\documentclass[12pt,a4paper]{article}
\usepackage{amssymb,amsmath}
\usepackage{amsfonts}
\usepackage{amsthm}
\usepackage{a4}
\usepackage{tikz}
\usepackage{verbatim}
\usepackage{caption}
\usepackage{float}
\newtheorem{prop}{Proposition}
\newtheorem*{theorem*}{Theorem}
\newtheorem{theorem}{Theorem}
\newtheorem{lemma}{Lemma}
\newtheorem{cor}{Corollary}
\newtheorem{remark}{Remark}
\newcommand{\F}{\mathbb{F}_{q}}
\newcommand{\ei}{\varepsilon_{\mathcal{C}}^{(i)}(H,K)}

\newcommand{\C}{\operatorname{Cen}}

\begin{document}
\baselineskip 18pt \title{\bf \ Semisimple finite group algebra of a generalized strongly monomial group}
\author{ Gurmeet K. Bakshi and Gurleen Kaur{\footnote {Research supported by NBHM(Ref No: 2/39(2)/2015/NBHM/R\&D-II/7424), Department of Atomic Energy, Government of India, is gratefully acknowledged.} \footnote{Corresponding author}} \\ {\em \small Centre for Advanced Study in
Mathematics,}\\
{\em \small Panjab University, Chandigarh 160014, India}\\{\em
\small email: gkbakshi@pu.ac.in and gurleenkaur992gk@gmail.com  } }
\date{}
{\maketitle}
\begin{abstract}  The complete algebraic structure of  semisimple finite group algebra of a generalized strongly monomial group is provided. This work extends the work of Broche and del R{\'{\i}}o on strongly monomial groups. The theory is complimented by an algorithm and is illustrated with an example.   \end{abstract}
 \vspace{.25cm}   \section{Introduction}

Let $\F$ be a  finite field with $q$ elements and $G$ a finite group of order coprime to  $q$. Let $\F G$ be the group algebra of $G$ over $\F$. During the last few years, lot of interest has been seen in understanding the primitive central idempotents and the Wedderburn decomposition of finite semisimple group algebra $\F G$, thus enabling coding theorist to determine the ideals of $\mathbb{F}_{q}G$ which are precisely all the group codes (e.g., see \cite{BGP2}, \cite{BGP1}, \cite{BGP3}, \cite{BRm2}, \cite{BRS}, \cite{Bs}, \cite{BR}, \cite{BP}, \cite{FLP}, \cite{FM}, \cite{GM}, \cite{Mr}, \cite{SBDR2}, \cite{SBDR}). \par If $G$ is abelian, this problem has been dealt in series of papers by several authors. Moving further to non abelian groups, a major step was taken by Broche and del R{\'{\i}}o in \cite{BR} where they provided the description of primitive central idempotents and the Wedderburn decomposition of $\F G$ when $G$ is a strongly monomial group. Recall that all abelian-by-supersolvable groups are strongly monomial and all strongly monomial groups are monomial. \par  Recently, in \cite{BK2}, we have defined generalized strongly monomial groups and proved that, beside strongly monomial groups, the class of generalized strongly monomial groups contain several important families of groups such as subnormally monomial groups (in particular Frobenius monomial groups) and solvable groups with all Sylow subgroups abelian extended by supersolvable groups. In this paper, we extend the work of Broche and del R{\'{\i}}o \cite{BR} to generalized strongly monomial groups. In section 2, we show (Theorem \ref{t1}) that if $(H,K)$ is a generalized strong Shoda pair of $G$ and $\mathcal{C}$ is a $q$-cyclotomic class of $\operatorname{Irr}(H/K)$ which contain generators of $\operatorname{Irr}(H/K)$, then the pair $((H, K), \mathcal{C})$ correspond to a primitive central idempotent of $\mathbb{F}_q G$. In this case, we also describe the structure of the corresponding simple component of $\F G$. This result allows us to describe the complete set of primitive central idempotents and the precise Wedderburn decomposition of semisimple finite group algebra of a generalized strongly monomial group (Corollary \ref{c2}). In section 3, we write a precise algorithm using the theory developed in section 2. In section 4, we illustrate the theory by computing the Wedderburn decomposition of semisimple group algebra $\mathbb{F}_{q}G$ for a generalized strongly monomial group $G$ which is not strongly monomial.
\section{Wedderburn decomposition} We begin by recalling the basic terminology related to Shoda pairs, strong Shoda pairs and generalized strong Shoda pairs. A pair $(H,K)$ of subgroups of $G$ is called a \textit{Shoda pair} of $G$ (see \cite{JdR}, Chapter 3) if the following conditions hold:\begin{description}\item[(i)] $K\unlhd H$, $H/K$ is cyclic;\item[(ii)] if $g \in G$ and $[H,g]\cap H\subseteq K$, then $g \in H$.\end{description}
For a Shoda pair $(H,K)$ of $G$, recall the following standard notations:  $$\widehat{H}:=\frac{1}{|H|}\displaystyle\sum_{h \in H}h,$$ $$\varepsilon(H,K):=\left\{\begin{array}{ll}\widehat{K}, & \hbox{$H=K$;} \\\prod(\widehat{K}-\widehat{L}), & \hbox{otherwise,}\end{array}\right.$$ where $L$ runs over the normal subgroups of $H$ minimal with respect to the property of including $K$ properly, and $$e(G,H,K):= {\rm~the~sum~of~all~the~distinct~}G{\rm {\operatorname{-}} conjugates~of~}\varepsilon(H,K).$$ \par
A Shoda pair $(H,K)$ of $G$ is called a \textit{strong Shoda pair} of $G$ (see \cite{JdR}, Chapter 3) if the following conditions hold:\begin{description}\item [(i)]$H \unlhd \operatorname{Cen}_{G}(\varepsilon(H,K))$;\item [(ii)] the distinct $G$-conjugates of $\varepsilon(H,K)$ are mutually orthogonal.\end{description} \par
More generally, a Shoda pair $(H,K)$ of $G$ is called a generalized strong Shoda pair of $G$ (see \cite{BK2}), if there is a chain $H=H_{0}\leq H_{1}\leq \cdots \leq H_{n}=G$ of subgroups of $G$ (called strong inductive chain from $H$ to $G$) such that the following conditions hold for all $0 \leq i < n$:  \begin{description} \item [(i)] $H_{i} \unlhd \operatorname{Cen}_{H_{i+1}}(\varepsilon^{(i)}(H,K))$;
\item [(ii)] the distinct $H_{i+1}$-conjugates of $\varepsilon^{(i)}(H,K)$ are mutually orthogonal, \end{description} where  $\varepsilon^{(0)}(H,K)= \varepsilon(H,K)$ and $\varepsilon^{(i)}(H,K)$ is the sum of all the distinct $H_{i}$-conjugates of $\varepsilon^{(i-1)}(H,K)$ for $1 \leq i \leq n$.  For notational convenience, we will denote $\varepsilon^{(n)}(H,K)$ by $\mathfrak{e}(G,H,K)$. \par Let $(H,K)$ be a Shoda pair of $G$ and let $H=H_{0}\leq H_{1}\leq \cdots \leq H_{n}=G$ be subgroups of $G$. By the repeated application of Lemma 3 of \cite{BK2}, $\varepsilon^{(i)}(H,K)$ is a rational multiple of $e_{\mathbb{Q}}(\lambda^{H_{i}})$, where $\lambda$ is any complex linear character on $H$ with kernel $K$ and $e_{\mathbb{Q}}(\lambda^{H_{i}})$ is the primitive central idempotent of the rational group algebra $\mathbb{Q}H_{i}$ associated with the complex irreducible character $\lambda^{H_{i}}$. Therefore, the centralizers of $\varepsilon^{(i)}(H,K)$ and $e_{\mathbb{Q}}(\lambda^{H_{i}})$ in $H_{i+1}$ coincide. Hence the above definition of a generalized strong Shoda pair of $G$ is equivalent to that given in \cite{BK2}.\par Recall that every Shoda pair $(H,K)$ of $G$ realizes a primitive central idempotent of $\mathbb{Q}G$ namely  $e_{\mathbb{Q}}(\lambda^{G})$, where $\lambda$ is any complex linear character on $H$ with kernel $K$. Two Shoda pairs of $G$ are called equivalent if they realize the same primitive central idempotent of $\mathbb{Q}G$. A group $G$ is called strongly monomial if every primitive central idempotent of $\mathbb{Q}G$ is realized by a strong Shoda pair of $G$ and a group $G$ is called generalized strongly monomial if every primitive central idempotent of $\mathbb{Q}G$ is realized by a generalized strong Shoda pair of $G$. Two generalized strong Shoda pairs of $G$ are called equivalent if they are equivalent as Shoda pairs of $G$. A set of representatives of distinct equivalence classes of generalized strong Shoda pairs of $G$ is called a complete and irredundant set of generalized strong Shoda pairs of $G$. For details on strongly monomial groups  and generalized strongly monomial groups, we refer to \cite{BK2} and \cite{OdRS04}.\par Given a Shoda pair $(H,K)$ of $G$, $\operatorname{Irr}(H/K)$ denotes the set of irreducible characters on $H/K$ over $\overline{\mathbb{F}}_{q}$, the algebraic closure of $\mathbb{F}_{q}$. If $k=[H:K]$ and $\xi_{k}$ is a primitive $k$th root of unity in $\overline{\mathbb{F}}_{q}$, then there is a natural action of $\operatorname{Gal}(\mathbb{F}_{q}(\xi_{k})/\mathbb{F}_{q})$ on $\operatorname{Irr}(H/K)$ by composition. The orbits of $\operatorname{Irr}(H/K)$ under this action are called $q$-cyclotomic classes. Denote by $\mathcal{C}_{q}(H/K)$ those $q$-cyclotomic classes which contain generators of $\operatorname{Irr}(H/K)$. If $\mathcal{C} \in \mathcal{C}_{q}(H/K)$ and $\chi \in \mathcal{C}$, set: $$\varepsilon_{\mathcal{C}}(H,K):= |H|^{-1}\sum_{h \in H} \operatorname{tr}(\chi(hK))h^{-1}=[H:K]^{-1}\widehat{K}\sum_{X \in H/K}\operatorname{tr}(\chi(X))h_{X}^{-1},$$ where  $\operatorname{tr}=\operatorname{tr}(\mathbb{F}_{q}(\xi_{[H:K]})/\mathbb{F}_{q}$) is the trace of the field extension $\mathbb{F}_{q}(\xi_{[H:K]})/\mathbb{F}_{q}$ and $h_{X}$ is a representative of $X \in H/K$.\vspace{.2cm}\\
 Define $$e_{\mathcal{C}}(G,H,K):={\rm~the~sum~of~all~the~distinct~}G{\rm {\operatorname{-}} conjugates~of~}\varepsilon_{\mathcal{C}}(H,K).$$ \par If $\psi$ is a character of a subgroup $H$ of $G$ and $x \in G$, then $\psi^{x}$ is a character of $H^{x}=x^{-1}Hx$ given by $\psi^{x}(g)=\psi(xgx^{-1})$, $g \in H^{x}$. Denote by $\psi^{G}$, the character $\psi$ induced to $G$. For $\alpha \in \mathbb{Q}G$ and $x \in G,$ $\alpha^{x}=x^{-1}\alpha x$. \vspace{.2cm}\\We recall the following theorem due to Broche and del R{\'{\i}}o \cite{BR}:
\begin{theorem*}{\rm(\cite{BR}, Theorem 7)} Let $G$ be a finite group and $\mathbb{F}_{q}$ be a finite field with $q$ elements such that $\mathbb{F}_{q}G$ is semisimple. \begin{description}\item [1] Let $(H,K)$ be a strong Shoda pair of $G$ and $\mathcal{C} \in \mathcal{C}_{q}(H/K)$. Then $e_{\mathcal{C}}(G,H,K)$ is a primitive central idempotent of $\mathbb{F}_{q}G$ and $$\mathbb{F}_{q}Ge_{\mathcal{C}}(G,H,K)\cong M_{[G:H]}(\mathbb{F}_{q^{(o/[E:H])}}),$$ where $E=E_{G}(H/K)$ is the stabilizer of any element of $\mathcal{C}_{q}(H/K)$ under the natural action of $N_{G}(K)$ and $o$ is the multiplicative order of $q$ module $[H:K]$. \item [2] Let $X$ be a set of strong Shoda pairs of $G$. If every primitive central idempotent of $\mathbb{Q}G$, the rational group algebra,  is of the form $e(G,H,K)$ for $(H,K) \in X$ then every primitive central idempotent of $\mathbb{F}_{q}G$ is of the form $e_{\mathcal{C}}(G,H,K)$ for $(H,K) \in X$ and $\mathcal{C} \in \mathcal{C}_{q}(H/K)$. \end{description} \qed \end{theorem*} \noindent The above theorem raises the following natural question: \begin{quote}  {\it  If  $(H,K)$ is a Shoda pair of $G$ and $\mathcal{C} \in \mathcal{C}_{q}(H/K)$, then, does the pair $((H, K), \mathcal{C})$ correspond to a primitive central idempotent of $\mathbb{F}_{q}G$?  If yes, can we describe the corresponding simple component? } \end{quote} The first part of the above question is not very difficult to see. Given such a pair $((H, K), \mathcal{C})$, we can see  analogous to Theorem 2.1 of \cite{OdRS04} that there is a unique $\alpha \in \mathbb{F}_{q}$ such that $e_{\mathbb{F}_q}(\lambda^G)= \alpha e_{\mathcal{C}}(G,H,K)$ for an arbitrary $\lambda \in \mathcal{C}$. Thus there does exist a primitive central idempotent, namely $\alpha e_{\mathcal{C}}(G,H,K)$ of $\mathbb{F}_{q}G$ that is associated with the pair $((H, K), \mathcal{C})$. To answer the above question, the main problem lies in describing the simple component $\mathbb{F}_{q}Ge_{\mathcal{C}}(G,H,K)$. We'll show, in Theorem \ref{t1}, that the structure of this simple component can be described by proceeding in steps, when $(H,K)$ is a generalized strong Shoda pair of $G$. Let us  elaborate on how to proceed by steps. If $(H,K)$ is a generalized strong Shoda pair of $G$, $H=H_{0}\leq H_{1}\leq \cdots \leq H_{n}=G$ is a strong inductive chain from $H$ to $G$ and $\mathcal{C}\in \mathcal{C}_{q}(H/K)$, set $$\varepsilon_{\mathcal{C}}^{(0)}(H,K) = \varepsilon_{\mathcal{C}}(H,K),$$ $$\varepsilon_{\mathcal{C}}^{(1)}(H,K)={\rm~the~sum~of~all~the~distinct~}H_{1}{\rm {\operatorname{-}} conjugates~of~}\varepsilon_{\mathcal{C}}^{(0)}(H,K), $$  $$\varepsilon_{\mathcal{C}}^{(2)}(H,K)  ={\rm~the~sum~of~all~the~distinct~}H_{2}{\rm {\operatorname{-}} conjugates~of~}\varepsilon_{\mathcal{C}}^{(1)}(H,K),$$ and finally  $$\varepsilon_{\mathcal{C}}^{(n)}(H,K)  ={\rm~the~sum~of~all~the~distinct~}H_{n}{\rm {\operatorname{-}} conjugates~of~}\varepsilon_{\mathcal{C}}^{(n-1)}(H,K).$$  Denote the final step $\varepsilon_{\mathcal{C}}^{(n)}(H,K) $ by $\mathfrak{e}_{\mathcal{C}}(G,H,K)$.
Apparently, one thinks that the definition of $\mathfrak{e}_{\mathcal{C}}(G,H,K)$ depends on strong inductive chain from $H$ to $G$. However, we will see that this is not the case, i.e., the final step $\mathfrak{e}_{\mathcal{C}}(G,H,K)$ is invariant of strong inductive chain from $H$ to $G$ (see Remark \ref{r1}). Observe that if $(H,K)$ is a strong Shoda pair of $G$, then $\mathfrak{e}_{\mathcal{C}}(G,H,K)= {e}_{\mathcal{C}}(G,H,K)$ on taking strong inductive chain: $H \leq G$.
	\begin{theorem}{\label{t1}} Let $\mathbb{F}_{q}$ be a finite field of order $q$ and  $G$ a finite group of order coprime to $q$.  Let $(H,K)$ be a generalized strong Shoda pair of $G$ and $H=H_{0}\leq H_{1}\leq \cdots \leq H_{n}=G$ a strong inductive chain from $H$ to $G$. If $\mathcal{C} \in \mathcal{C}_{q}(H/K)$, then the following statements hold:
		\begin{description} \item [(i)] $\mathfrak{e}_{\mathcal{C}}(G,H,K)$ is a primitive central idempotent of $\mathbb{F}_{q} G$. More generally, for any $\lambda \in \mathcal{C}$, $\mathfrak{e}_{\mathcal{C}}(G,H,K)= e_{\mathbb{F}_{q}}(\lambda^{G})$; \item [(ii)]  $\mathbb{F}_{q}G \mathfrak{e}_{\mathcal{C}}(G,H,K)  \cong M_{[G:H]}(\mathbb{F}_{{q}^{l_{(H,K)}}}),$ where $l_{(H,K)}=\frac{o}{\prod_{0\leq i< n}[C_{i}:H_{i}]}$, $o$ is the multiplicative order of $q$ module $[H:K]$ and $C_{i} = \operatorname{Cen}_{H_{i+1}}(\varepsilon_{\mathcal{C}}^{(i)}(H,K))$.  \end{description}
	\end{theorem}
\begin{remark}\label{r1} Since $ \mathfrak{e}_{\mathcal{C}}(G,H,K)= e_{\mathbb{F}_{q}}(\lambda^{G})$, it follows immediately that the construction of $\mathfrak{e}_{\mathcal{C}}(G,H,K)$ is independent of strong inductive chain from $H$ to $G$.
\end{remark}\begin{remark}\label{r2} $\prod_{0\leq i  < n}[C_{i}: H_{i}]$ always divides $o$.\end{remark}
\begin{lemma}{\label{l1}}If $(H, K)$ is a Shoda pair of a finite group $G$ and $\lambda$ is an $\overline{\mathbb{F}}_{q}$-irreducible character on $H$ with kernel $K$, then $\lambda^G$ is irreducible.
	
\end{lemma}
\noindent{\bf Proof.} This follows from Corollary 45.4 of \cite{CR}. $\Box$
\begin{lemma}{\label{l2}} Let $S$ be a subgroup of a finite group $G$ and $\psi$ an $\overline{\mathbb{F}}_{q}$-irreducible character on $S$ such that $\psi^{G}$ is irreducible, then there exists $ \alpha \in \mathbb{F}_{q}$ such that $$e_{\mathbb{F}_{q}}(\psi^{G})=\alpha \sum_{g \in T}e_{\mathbb{F}_{q}}(\psi)^{g},$$ where $T$ is a right transversal of $\operatorname{Cen}_{G}(e_{\mathbb{F}_{q}}(\psi))$ in $G$. Furthermore $\alpha =1$, if the distinct $G$-conjugates of $e_{\mathbb{F}_{q}}(\psi)$ are mutually orthogonal and in this case, $\F G e_{\mathbb{F}_{q}}(\psi^{G})$ is isomorphic to $ M_{[G:C]}(\F C e_{\mathbb{F}_{q}}(\psi) )$, where
	$C= \C_{G}(e_{\mathbb{F}_{q}}(\psi))$.
\end{lemma} \noindent{\bf Proof.} The proof is analogous to that of (\cite{BK2}, Lemma 3, Proposition 2(i)).\qed\vspace{.2cm} \\ Throughout the rest of the paper, $(H,K)$  is  a generalized strong Shoda pair of $G$ and $H=H_{0}\leq H_{1}\leq \cdots \leq H_{n}=G$ is a strong inductive chain from $H$ to $G$.\par Let $p$ be the prime divisor of $q$ and let $\mathbb{Z}_{(p)}$ denote the localization of $\mathbb{Z}$ at $p$. Then $\mathbb{F}_{p}$ can be identified with the residue field of $\mathbb{Z}_{(p)}$. Consider the natural ring homomorphism from $\mathbb{Z}_{(p)}$ onto $\mathbb{F}_{p}$ and extend it linearly from $\mathbb{Z}_{(p)}G$ onto $\mathbb{F}_{p}G.$ Denote by $\overline{\alpha}$, the image of $\alpha \in \mathbb{Z}_{(p)}G$ under the above epimorphism. By (\cite{BK2}, Lemma 3), it follows that $\varepsilon^{(i)}(H,K)$ is a primitive central idempotent of $\mathbb{Z}_{(p)}H_{i}$ and hence $\overline{\varepsilon^{(i)}(H,K)}$ is a central idempotent of $\mathbb{F}_{p}H_{i}$ for all $0\leq i\leq n$.
\begin{lemma}{\label{l3}} If, for some $i$, $\overline{\varepsilon^{(i)}(H,K)}= e_{1} + \cdots + e_{l}$ is the decomposition of $\overline{\varepsilon^{(i)}(H,K)}$ into sum of  primitive central idempotents of $\mathbb{F}_{q}H_{i}$, then the following statements hold for all $ 1 \leq k \leq l$: \begin{description}   \item [(i)] $\operatorname{Cen}_{H_{i+1}}(e_{k})$ is a subgroup of $\operatorname{Cen}_{H_{i+1}}(\varepsilon^{(i)}(H,K))$; \item [(ii)] $H_{i} \unlhd  \operatorname{Cen}_{H_{i+1}}(e_{k})$;   \item [(iii)] the distinct $H_{i+1}$-conjugates of $e_{k}$ are mutually orthogonal.  \end{description}	
\end{lemma}\noindent{\bf Proof.} By renaming $e_{i}$'s, we may assume that $k=1$.  On multiplying the equation $\overline{\varepsilon^{(i)}(H,K)}= e_{1} + \cdots + e_{l}$ with $e_{1}$, we obtain that \begin{equation}{\label{e1} }
\overline{\varepsilon^{(i)}(H,K)}e_{1}= me_{1}= e_{1}\overline{\varepsilon^{(i)}(H,K)}, \end{equation} where $m$ is the number of $e_{j}$'s that are equal to $e_{1}$. Note that the characteristic of $\mathbb{F}_{q}$ does not divide $m$. Consider any $ g \in H_{i+1}\setminus \C_{H_{i+1}}(\varepsilon^{(i)}(H,K)) $. Conjugating  eqn (\ref{e1}) by $g$, we have
\begin{equation}{\label{e2}}
\overline{\varepsilon^{(i)}(H,K)^{g}} e_{1}^{g}= me_{1}^{g}= e_{1}^{g}\overline{\varepsilon^{(i)}(H,K)^{g}}. \end{equation} Since $\varepsilon^{(i)}(H,K) \varepsilon^{(i)}(H,K)^{g} =0$, we get from eqns (\ref{e1}) and (\ref{e2}) that $e_{1}e_{1}^{g}=0$. This gives that $\C_{H_{i+1}}(e_{1})$ is a subgroup of $\C_{H_{i+1}}(\varepsilon^{(i)}(H,K))$, which proves (i). The normality of $H_{i}$ in $\C_{H_{i+1}}(\varepsilon^{(i)}(H,K))$ immediately yields (ii), in view of (i). Since $\operatorname{Cen}_{H_{i+1}}(e_{1})\leq \operatorname{Cen}_{H_{i+1}}(\varepsilon^{(i)}(H,K))$ and $e_{1}e_{1}^{g}=0$ if $g \notin \operatorname{Cen}_{H_{i+1}}(\varepsilon^{(i)}(H,K))$, to prove (iii), we only  need to show that $e_{1} e_{1}^{g} =0$  if $ g \in \C_{H_{i+1}}(\varepsilon^{(i)}(H,K)) \setminus \C_{H_{i+1}}(e_{1}) $.  But this clearly holds because if $ g \in  \C_{H_{i+1}}(\varepsilon^{(i)}(H,K)) \setminus \C_{H_{i+1}}(e_1),$ then $e_{1}$ and $e_{1}^{g}$ are distinct primitive central idempotents of $\F H_{i}$, as $H_{i}\unlhd \operatorname{Cen}_{H_{i+1}}(\varepsilon^{(i)}(H,K))$. Hence, the lemma is proved. \qed
\begin{prop}{\label{p1}}  There exist subsets $R_i$,  $ 0 \leq i \leq n$, of $ \mathcal{C}_{q}(H/K)$ with \break $R_{0} = \mathcal{C}_{q}(H/K)$ and $R_{i}$ contained in $R_{i-1}$ such that $\operatorname{Cen}_{H_{i}}(\varepsilon^{(i-1)}(H,K))$ acts on $\{\varepsilon_{\mathcal{C}}^{(i-1)}(H,K)~|~\mathcal{C} \in R_{i-1}\}$ by conjugation and $\{\varepsilon_{\mathcal{C}}^{(i-1)}(H,K)~|~\mathcal{C} \in R_{i}\}$ is a set of representatives under this action for all $1 \leq i \leq n$. Furthermore, as a consequence, the following statements hold for all $0\leq i\leq n$:
	\begin{description}\item [(i)] $\overline{\varepsilon^{(i)}(H,K)} = \sum_{\mathcal{C} \in R_{i}}\varepsilon_{\mathcal{C}}^{(i)}(H,K)$;\item [(ii)] $\varepsilon_{\mathcal{C}}^{(i)}(H,K)$, $ \mathcal{C} \in R_{i}$, are distinct primitive central idempotents of $\mathbb{F}_{q}H_{i}$.\break Moreover, $\ei= e_{\F}(\lambda^{H_{i}}) $ for an arbitrary $\lambda \in \mathcal{C}$.
	\end{description}
\end{prop}\noindent{\bf Proof.}   We will prove the result by induction on $i$. Suppose $i=0$. Then $\varepsilon_{\mathcal{C}}^{(0)}(H,K)= \varepsilon_{\mathcal{C}}(H,K)$ and the result follows from (\cite{BR}, Corollary 3, Lemma 6). Assume that the result is true for $i=k$, where $0 \leq k < n$.  We'll prove the result for $i= k+1$. Let $ g \in \operatorname{Cen}_{H_{k+1}}(\varepsilon^{(k)}(H,K))$. Then
$$\overline{\varepsilon^{(k)}(H,K)} = \sum_{\mathcal{C} \in R_{k}}\varepsilon_{\mathcal{C}}^{(k)}(H,K)=\sum_{\mathcal{C} \in R_k}(\varepsilon_{\mathcal{C}}^{(k)}(H,K))^{g} .$$
By induction, $\varepsilon_{\mathcal{C}}^{(k)}(H,K)$, $ \mathcal{C} \in R_{k}$, are distinct primitive central idempotents of $\mathbb{F}_{q}H_{k}$. Since $ H_{k} \unlhd \operatorname{Cen}_{H_{k+1}}(\varepsilon^{(k)}(H,K))$, we also have that  $\varepsilon_{\mathcal{C}}^{(k)}(H,K)^{g}$, $ \mathcal{C} \in R_k$,  are also distinct primitive central idempotents of $\mathbb{F}_{q}H_{k}$. But  a central idempotent can be uniquely written as a sum of primitive central idempotents, therefore $$ \{ (\varepsilon_{\mathcal{C}}^{(k)}(H,K))^{g} ~| ~ \mathcal{C} \in R_{k} \} =   \{ \varepsilon_{\mathcal{C}}^{(k)}(H,K)~ | ~ \mathcal{C} \in R_{k} \} .$$ This gives an action of $\operatorname{Cen}_{H_{k+1}}(\varepsilon^{(k)}(H,K))$ on $\{ \varepsilon_{\mathcal{C}}^{(k)}(H,K)~|~  \mathcal{C} \in R_k\}$ by \break conjugation. Under this action, let $R_{k+1}$ contained in $R_{k}$ be such that \break $\{ \varepsilon_{\mathcal{C}}^{(k)}(H,K)~ | ~ \mathcal{C} \in R_{k+1} \} $ is a set of representatives of distinct orbits. For any $ \mathcal{C} \in R_{k+1}$, let $E_{\mathcal{C}}$ be the stabilizer of $\varepsilon_{\mathcal{C}}^{(k)}(H,K)$ in $ \operatorname{Cen}_{H_{k+1}}(\varepsilon^{(k)}(H,K))$. In view of Lemma \ref{l3},  we have $E_{\mathcal{C}} = \operatorname{Cen}_{H_{k+1}}(\varepsilon_{\mathcal{C}}^{(k)}(H,K))$. Let
$T_{\mathcal{C}}$ be a right transversal of $E_{\mathcal{C}}$ in $\operatorname{Cen}_{H_{k+1}}(\varepsilon^{(k)}(H,K))$. Let $T$ be a right transversal of
  $\operatorname{Cen}_{H_{k+1}}(\varepsilon^{(k)}(H,K))$ in $H_{k+1}$  so that $\{xy~|~x \in T_{\mathcal{C}}, y \in T\}$ is a right transversal of $E_{\mathcal{C}}$ in $H_{k+1}$. Now,
\begin{align*}
 \overline{\varepsilon^{(k+1)}(H,K)} &= \sum_{y \in T}(\overline{\varepsilon^{(k)}(H,K)})^{y} \\
 &= \sum_{y \in T} \sum_{\mathcal{C} \in R_{k}}(\varepsilon_{\mathcal{C}}^{(k)}(H,K))^{y}~{\rm (by~induction)}\\ &=
 \sum_{y \in T}( \sum_{\mathcal{C} \in R_{k+1}}\sum_{x \in T_{\mathcal{C}}}(\varepsilon_{\mathcal{C}}^{(k)}(H,K))^x)^y \\&= \sum_{\mathcal{C} \in R_{k+1}} \sum_{y \in T}\sum_{x \in T_{\mathcal{C}}}(\varepsilon_{\mathcal{C}}^{(k)}(H,K))^{xy} \\  &=\sum_{\mathcal{C} \in R_{k+1}}\varepsilon_{\mathcal{C}}^{(k+1)}(H,K).
\end{align*}
 \vspace{.2cm} \\
This proves (i)  for $i=k+1$. Next, we prove (ii) for $ i= k+1$. Consider any $ \mathcal{C} \in R_{k+1}$ and an arbitrary $\lambda \in \mathcal{C} $. By  Lemma \ref{l1}, $\lambda^{H_{k+1}}$ is irreducible and by Lemma \ref{l2}, $$e_{\F}(\lambda^{H_{k+1}}) = \alpha ({\rm the~sum~of~all~the~distinct~}H_{k+1}{\rm {\operatorname{-}} conjugates~of~} e_{\F}(\lambda^{H_{k}})),$$ for some $\alpha \in \F$. By induction, $e_{\F}(\lambda^{H_{k}}) =  \varepsilon_{\mathcal{C}}^{(k)}(H,K)$. Also, by induction, both (i) and (ii) hold for $i=k$, so we obtain using Lemma \ref{l3} that the distinct $H_{k+1}$-conjugates of $\varepsilon_{\mathcal{C}}^{(k)}(H,K)$ are mutually orthogonal and thus, by Lemma \ref{l2}, $\alpha =1$. This gives that $e_{\F}(\lambda^{H_{k+1}})$ equals ${\rm the~sum~of~all~the~distinct~}H_{k+1}{\rm {\operatorname{-}} conjugates~of~}$  $ \varepsilon_{\mathcal{C}}^{(k)}(H,K)$, i.e., $e_{\F}(\lambda^{H_{k+1}})= \varepsilon_{\mathcal{C}}^{(k+1)}(H,K).$ Thus $\varepsilon_{\mathcal{C}}^{(k+1)}(H,K)$ is a primitive central idempotent of $\F H_{k+1}$. Next, we show that if $\mathcal{C}' \in R_{k+1}$ is different from $\mathcal{C}$, then $\varepsilon_{\mathcal{C'}}^{(k+1)}(H,K)$ and $\varepsilon_{\mathcal{C}}^{(k+1)}(H,K)$ are mutually orthogonal. Let $m$ be the number of $\mathcal{C}'$ in $R_{k+1}$ such that $\varepsilon_{\mathcal{C}'}^{(k+1)}(H,K) = \varepsilon_{\mathcal{C}}^{(k+1)}(H,K)$. On multiplying $\overline{\varepsilon^{(k+1)}(H,K)}=\sum_{\mathcal{C}\in R_{k+1}}\varepsilon_{\mathcal{C}}^{(k+1)}(H,K)$ with $\varepsilon_{\mathcal{C}}^{(k+1)}(H,K)$, we get that $\overline{\varepsilon^{(k+1)}(H,K)}\varepsilon_{\mathcal{C}}^{(k+1)}(H,K) = m \varepsilon_{\mathcal{C}}^{(k+1)}(H,K)$. Since both  $\overline{\varepsilon^{(k+1)}(H,K)}\varepsilon_{\mathcal{C}}^{(k+1)}(H,K) $ and $\varepsilon_{\mathcal{C}}^{(k+1)}(H,K)$ are idempotents, it follows that $m=1$. This proves (ii) for $i=k+1$ and hence completes the proof.  \qed \begin{cor}{\label{c1}} For every $ \mathcal{C} \in  \mathcal{C}_{q}(H/K)$ and for every $i$, $0\leq i\leq n$, there exists $\mathcal{C}_{i} \in R_{i} $ such that $ \ei = \varepsilon_{\mathcal{C}_{i}}^{(i)}(H,K)$ and consequently, the following properties hold: \begin{description}\item [(i)] $H_{i}\unlhd \operatorname{Cen}_{H_{i+1}}(\varepsilon_{\mathcal{C}}^{(i)}(H,K)),~i\neq n$;\item [(ii)] the distinct $H_{i+1}$-conjugates of $\varepsilon_{\mathcal{C}}^{(i)}(H,K)$ are mutually orthogonal,~$i\neq n.$\end{description} Moreover, $\varepsilon_{\mathcal{C}}^{(i)}(H,K)=e_{\mathbb{F}_{q}}(\lambda^{H_{i}})$ for an arbitrary $\lambda \in \mathcal{C}$.\end{cor}\noindent{\bf Proof.} We assert by induction on $i$ that for every $ \mathcal{C} \in  \mathcal{C}_{q}(H/K)$, there exists $\mathcal{C}_{i} \in R_{i} $ such that $ \ei = \varepsilon_{\mathcal{C}_{i}}^{(i)}(H,K)$. The result trivially holds true for $i=0$, as $R_{0}= \mathcal{C}_{q}(H/K)$. Suppose that the result is true for $i=k$. Hence there is a $C_{k} \in R_{k}$ such that $ \varepsilon_{\mathcal{C}}^{(k)}(H,K) = \varepsilon_{\mathcal{C}_k}^{(k)}(H,K)$. Thus the sum of all $H_{k+1}$-conjugates of $\varepsilon_{\mathcal{C}}^{(k)}(H,K)$ and $\varepsilon_{\mathcal{C}_{k}}^{(k)}(H,K)$ are also same, i.e., $ \varepsilon_{\mathcal{C}}^{(k+1)}(H,K) = \varepsilon_{\mathcal{C}_{k}}^{(k+1)}(H,K)$. We now show that there is a $C_{k+1} \in R_{k+1}$ such that $ \varepsilon_{\mathcal{C}_{k}}^{(k+1)}(H,K) =  \varepsilon_{\mathcal{C}_{k+1}}^{(k+1)}(H,K)$.  Recall from Proposition \ref{p1} that  $\{ \varepsilon_{\mathcal{C}}^{(k)}(H,K)~ | ~ \mathcal{C} \in R_{k+1} \} $ is a set of representatives of distinct orbits under the  action of $\operatorname{Cen}_{H_{k+1}}(\varepsilon^{(k)}(H,K))$ on $\{ \varepsilon_{\mathcal{C}}^{(k)}(H,K)~|~  \mathcal{C} \in R_{k}\}$ by conjugation.  So there is a $C_{k+1} \in R_{k+1}$ such that  $\varepsilon_{\mathcal{C}_{k}}^{(k)}(H,K) = (\varepsilon_{\mathcal{C}_{k+1}}^{(k)}(H,K))^g$ for some $ g \in \operatorname{Cen}_{H_{k+1}}(\varepsilon^{(k)}(H,K))$. Thus the sum of all $H_{k+1}$-conjugates of $\varepsilon_{\mathcal{C}_{k}}^{(k)}(H,K)$ is equal to the sum of all $H_{k+1}$-conjugates of $(\varepsilon_{\mathcal{C}_{k+1}}^{(k)}(H,K))^g$, which is same as the sum of all $H_{k+1}$-conjugates of $\varepsilon_{\mathcal{C}_{k+1}}^{(k)}(H,K)$,  as $ g \in H_{k+1}$. Hence  $ \varepsilon_{\mathcal{C}_k}^{(k+1)}(H,K) =  \varepsilon_{\mathcal{C}_{k+1}}^{(k+1)}(H,K)$, as desired. This proves the assertion. Furthermore, (i) and (ii) hold using Proposition 1 and Lemma 3. Moreover, in view of (ii), the repeated application of Lemma \ref{l2} yields that $\varepsilon_{\mathcal{C}}^{(i)}(H,K)=e_{\mathbb{F}_{q}}(\lambda^{H_{i}})$ for an arbitrary $\lambda \in \mathcal{C}$ and $0\leq i\leq n$.\qed \vspace{.5cm}\\Denote by $R*_{\tau}^{\sigma}G$, the crossed product of the group $G$ over the ring $R$ with action $\sigma$ and twisting $\tau$ ( for details on crossed product, see (\cite{JdR}, Chapter 2)).
\begin{prop}{\label{p2}} The following statements hold for all $\mathcal{C}\in \mathcal{C}_{q}(H/K)$ and \linebreak $0\leq i < n$: \begin{description} 	 \item [(i)] $\mathbb{F}_{q}H_{i+1}\varepsilon_{\mathcal{C}}^{(i+1)}(H,K) \cong M_{[H_{i+1}:C_{i}]} (\mathbb{F}_{q} C_{i} \varepsilon_{\mathcal{C}}^{(i)}(H,K))$, where $ C_{i} = \operatorname{Cen}_{H_{i+1}}(\varepsilon^{(i)}_{\mathcal{C}}(H,K)).$	\item [(ii)] $\mathbb{F}_{q} C_{i} \varepsilon_{\mathcal{C}}^{(i)}(H,K)$ $\cong$ $\mathbb{F}_{q} H_{i} \varepsilon_{\mathcal{C}}^{(i)}(H,K) *^{\sigma_{i}}_{\tau_{i}} C_{i}/H_{i}$, where the action \break $\sigma_{i}:C_{i}/H_{i}\rightarrow \operatorname{Aut}(\mathbb{F}_{q}H_{i}\varepsilon_{\mathcal{C}}^{(i)}(H,K))$ maps $x$ to the conjugation automorphism $(\sigma_{i})_{x}$ on $\mathbb{F}_{q}H_{i}\varepsilon_{\mathcal{C}}^{(i)}(H,K)$ induced by the fixed inverse image $\overline{x}$ of $x$ under the natural map $C_{i}\rightarrow C_{i}/H_{i}$ and the twisting $\tau_{i}:$~$ C_{i}/H_{i}\times C_{i}/H_{i}\rightarrow \mathcal{U}(\mathbb{F}_{q}H_{i}\varepsilon_{\mathcal{C}}^{(i)}(H,K))$ is given by $\tau_{i}(x,y)=h\varepsilon_{\mathcal{C}}^{(i)}(H,K)$, where $h \in H_{i}$ is such that $\overline{x}.\overline{y}=h\overline{xy}$ and $\mathcal{U}(\mathbb{F}_{q}H_{i}\varepsilon_{\mathcal{C}}^{(i)}(H,K))$ is the unit group of $\mathbb{F}_{q}H_{i}\varepsilon_{\mathcal{C}}^{(i)}(H,K)$. \item [(iii)] $(\sigma_{i})_{x}$ is not an inner automorphism of $ \mathbb{F}_{q} H_{i} \varepsilon_{\mathcal{C}}^{(i)}(H,K)$ for every non identity $x \in C_{i}/H_{i}$.  \item [(iv)] $C_{i}/H_{i}$ acts faithfully on the center $\mathcal{Z}(\mathbb{F}_{q}H_{i} \varepsilon_{\mathcal{C}}^{(i)}(H,K))$ of $\mathbb{F}_{q} H_{i} \varepsilon_{\mathcal{C}}^{(i)}(H,K)$ by conjugation. \end{description} \end{prop}
\noindent {\bf Proof.} Let $\mathcal{C} \in \mathcal{C}_{q}(H/K)$ and $0 \leq i < n$. Consider an arbitrary $\lambda \in  \mathcal{C}$.  \\(i) In view of Corollary \ref{c1}, $e_{\F}(\lambda^{H_{i}})=\varepsilon_{\mathcal{C}}^{(i)}(H,K),$ $e_{\F}(\lambda^{H_{i+1}})=\varepsilon_{\mathcal{C}}^{(i+1)}(H,K)$ and the distinct $H_{i+1}$-conjugates of $e_{\F}(\lambda^{H_{i}})$ are mutually orthogonal. Hence, by taking $S=H_{i}$, $ G = H_{i+1}$ and $\psi = \lambda^{H_{i}}$ in Lemma \ref{l2}, we get $$\mathbb{F}_{q}H_{i+1}\varepsilon_{\mathcal{C}}^{(i+1)}(H,K) \cong M_{[H_{i+1}: C_{i}]}(\mathbb{F}_{q}C_{i}\varepsilon_{\mathcal{C}}^{(i)}(H,K)).$$  This proves (i).\\
(ii) In view of Corollary \ref{c1}, we have $H_i \unlhd C_i$. Therefore, (ii) follows. \\
(iii) Suppose  $ x \in C_{i}/H_{i}$ is such that $(\sigma_{i})_{x}$ is an inner automorphism of $ \mathbb{F}_{q}H_{i}\varepsilon_{\mathcal{C}}^{(i)}(H,K)$. So there is a unit $u$ in $ \mathbb{F}_{q}H_{i}\varepsilon_{\mathcal{C}}^{(i)}(H,K)$ such that \begin{equation}\label{e3}
\overline{x}\alpha \overline{x}^{-1}= u\alpha u^{-1}~~ \forall ~\alpha \in \mathbb{F}_{q} H_{i} \varepsilon_{\mathcal{C}}^{(i)}(H,K).
\end{equation}
In particular,
 \begin{equation}{\label{e4} }
\overline{x}h\ei \overline{x}^{-1}= u h\ei u^{-1}~~ \forall~  h \in H_{i}.
\end{equation} Let $\rho_i$ be the representation of $H_{i}$ affording the character $\lambda^{H_{i}}$. Extending $\rho_{i}$ linearly on $\mathbb{F}_{q}H_{i}$ and applying on eqn (\ref{e4}), we get
$$\rho_{i}(\overline{x}h\ei \overline{x}^{-1})= \rho_i(uh\ei u^{-1})~~ \forall~  h \in H_{i}.$$ Since $\rho_{i}$ is a ring homomorphism and $\rho_i(\ei)$ is the identity matrix, we obtain $$\rho_{i}(\overline{x}h \overline{x}^{-1})= \rho_i(u)\rho_{i}(h)\rho_{i}(u^{-1})~~ \forall~  h \in H_{i}.$$
On taking trace, $$(\lambda^{H_{i}})^{\overline{x}}(h) = \lambda^{H_{i}}(h)~~ \forall~  h \in H_{i}.$$ Thus $\overline{x}$ belongs to the inertia group of $\lambda^{H_{i}}$ in $C_{i}.$ Since $H_i \unlhd C_i$ and $\lambda^{C_{i}}$ is irreducible, it follows from Mackey's irreducibility criterion that $\overline{x} \in H_{i}$, i.e., $x$ is identity. \\
(iv) Clearly,  $C_{i}/H_{i}$ acts  on  $\mathcal{Z}(\mathbb{F}_{q} H_{i} \varepsilon_{\mathcal{C}}^{(i)}(H,K))$ by conjugation. All we need to see is that the action is faithful. Let $ x \in C_{i}/H_{i}$ be such that
\begin{equation}{\label{e5} }
\overline{x}\alpha \overline{x}^{-1}= \alpha ~~ \forall ~\alpha \in \mathcal{Z}(\mathbb{F}_{q} H_{i} \varepsilon_{\mathcal{C}}^{(i)}(H,K)).
\end{equation}
Consider an arbitrary $ h \in H_{i}$. Let $T$ be a right transversal of $\C_{H_{i}}(h)$ in $H_{i}$. Observe that $\sum_{ g \in T}h^{g}$ commutes with all the elements of $H_{i}$ and hence it belongs to
$\mathcal{Z}(\mathbb{F}_{q} H_{i} \varepsilon_{\mathcal{C}}^{(i)}(H,K))$. Thus, from eqn (\ref{e5}), $$\overline{x}(\sum_{ g \in T}h^g) \overline{x}^{-1}= \sum_{ g \in T}h^g .$$ Let $\rho_{i}$ be the representation afforded by $\lambda^{H_{i}}$. Extending $\rho_{i}$ linearly on $\mathbb{F}_{q}H_{i}$, we have $$\rho_{i}(\overline{x}(\sum_{ g \in T}h^g) \overline{x}^{-1})= \rho_{i}(\sum_{ g \in T}h^g).$$ On taking trace, it follows that $(\lambda^{H_{i}})^{\overline{x}}(h) = \lambda^{H_{i}}(h)$. This is true for all $h \in H_{i}$. Hence
$\overline{x}$ belongs to the inertia group of $\lambda^{H_{i}}$ in $C_{i} $, which gives $\overline{x} \in H_{i}$ and hence $x$ is identity.  \qed \vspace{.5cm}\\
\noindent{\bf Proof of Theorem \ref{t1}.} From Corollary \ref{c1}, $e_{\mathbb{F}_{q}}(\lambda^{H_{i}})=\varepsilon_{\mathcal{C}}^{(i)}(H,K)$ for $0 \leq i \leq n$. In particular, $i=n$ gives that $e_{\mathbb{F}_{q}}(\lambda^{G})$ equals $\varepsilon_{\mathcal{C}}^{(n)}(H,K)$, i.e., $\mathfrak{e}_{\mathcal{C}}(G,H,K)$. This proves (i). We now proceed to prove (ii). We will show that the following holds for all $ 1 \leq i \leq n$:
\begin{equation}{\label{e6}}\mathbb{F}_{q}H_{i}\varepsilon_{\mathcal{C}}^{(i)}(H,K) \cong M_{[H_{i}: H]}(\mathbb{F}_{q^{l_{i}}}),\end{equation} where $l_{i}= \frac{o}{[C_{0}: H_{0}][C_{1}: H_{1}] \cdots [C_{i-1}: H_{i-1}] }.$  Then (ii) follows from eqn (\ref{e6}) by taking $i=n$. To prove eqn (\ref{e6}), we will  induct on $i$.  For $i=1$, the result follows from  Theorem 7 of \cite{BR}, as $(H,K)$ is a strong Shoda pair of $H_{1}$ and $\varepsilon_{\mathcal{C}}^{(1)}(H,K) = e_{\mathcal{C}}(H_1,H,K)$.  Suppose that the result is true for $i=k$, where $1 \leq k < n$, i.e., \begin{equation}{\label{e7}}
\mathbb{F}_{q}H_k \varepsilon_{\mathcal{C}}^{(k)}(H,K) \cong M_{[H_k: H]}(\mathbb{F}_{q^{l_k}}).
\end{equation} Using Proposition  \ref{p2},
\begin{equation}{\label{e8}}
\mathbb{F}_{q}H_{k+1}  \varepsilon_{\mathcal{C}}^{(k+1)}(H,K)  \cong  M_{[H_{k+1}: C_{k}]}(\mathbb{F}_q H_k \varepsilon_{\mathcal{C}}^{(k)}(H,K) *^{\sigma_k}_{\tau_k} C_k/H_k). \end{equation} Since  $\mathbb{F}_{q}H_{k}\varepsilon_{\mathcal{C}}^{(k)}(H,K)$ is a simple ring and, by Proposition \ref{p2}, $(\sigma_{k})_{x}$ is not an inner automorphism of $ \mathbb{F}_{q}H_{k}\varepsilon_{\mathcal{C}}^{(k)}(H,K)$ for every non identity $x \in C_k/H_k$, it follows from Lemma 2.6.1 of \cite{JdR} that $\mathbb{F}_{q}H_{k}\varepsilon_{\mathcal{C}}^{(k)}(H,K)*^{\sigma_{k}}_{\tau_{k}}C_{k}/H_{k}$ is also a simple ring and its center  $\mathcal{Z}(\mathbb{F}_{q}H_{k}\varepsilon_{\mathcal{C}}^{(k)}(H,K) *^{\sigma_{k}}_{\tau_{k}}C_{k}/H_{k})$ equals $\mathcal{Z}(\mathbb{F}_{q}H_{k} \varepsilon_{\mathcal{C}}^{(k)}(H,K))^{C_{k}/H_{k}}$, the fixed subring of $\mathcal{Z}(\mathbb{F}_{q}H_{k}\varepsilon_{\mathcal{C}}^{(k)}(H,K))$. Consequently,  $\mathbb{F}_{q}H_{k}\varepsilon_{\mathcal{C}}^{(k)}(H,K)*^{\sigma_{k}}_{\tau_{k}}C_{k}/H_{k}$ is a matrix ring over $\mathcal{Z}(\mathbb{F}_{q}H_{k} \varepsilon_{\mathcal{C}}^{(k)}(H,K))^{C_{k}/H_{k}}$, say of size  $m \times m$. This gives
$$ \mathbb{F}_{q}H_{k+1}\varepsilon_{\mathcal{C}}^{(k+1)}(H,K) \cong  M_{m[H_{k+1}: C_{k}]}(\mathcal{Z}(\mathbb{F}_{q}H_{k} \varepsilon_{\mathcal{C}}^{(k)}(H,K))^{C_k/H_k}).$$  The faithful action of $C_k/H_k$ on
$\mathcal{Z}(\mathbb{F}_q H_k \varepsilon_{\mathcal{C}}^{(k)}(H,K))$ (Proposition \ref{p2}) gives that
$ \operatorname{dim}_{\mathbb{F}_q}(\mathcal{Z}(\mathbb{F}_q H_k \varepsilon_{\mathcal{C}}^{(k)}(H,K))^{C_k/H_k}) = \frac{\operatorname{dim}_{\mathbb{F}_q}(\mathcal{Z}(\mathbb{F}_q H_k \varepsilon_{\mathcal{C}}^{(k)}(H,K)))}{[C_k: H_k]}$.  From eqn (\ref{e7}), \break $\operatorname{dim}_{\mathbb{F}_q}(\mathcal{Z}(\mathbb{F}_q H_k \varepsilon_{\mathcal{C}}^{(k)}(H,K)))=  l_k$. Therefore, $ \operatorname{dim}_{\mathbb{F}_q}(\mathcal{Z}(\mathbb{F}_q H_k \varepsilon_{\mathcal{C}}^{(k)}(H,K))^{C_k/H_k}) = \frac{l_k}{[C_k: H_k]} =  l_{k+1}$.  This gives \begin{equation}{\label{e9}} \mathbb{F}_{q}H_{k+1}  \varepsilon_{\mathcal{C}}^{(k+1)}(H,K) \cong  M_{m[H_{k+1}: C_{k}]}(\mathbb{F}_{q^{l_{k+1}}}).  \end{equation}  On comparing the dimension of $\mathbb{F}_{q}H_{k+1}  \varepsilon_{\mathcal{C}}^{(k+1)}(H,K)$ over $ \mathbb{F}_{q}$ in eqns (\ref{e8}) and (\ref{e9}), we get $ [H_{k+1} : C_k]^2 [C_k : H_k]  \operatorname{dim}_{\mathbb{F}_q}(\mathbb{F}_q H_k \varepsilon_{\mathcal{C}}^{(k)}(H,K) ) = m^2 [H_{k+1} : C_k]^2 l_{k+1}$. Since from eqn (\ref{e7}),  $ \operatorname{dim}_{\mathbb{F}_q}(\mathbb{F}_q H_k \varepsilon_{\mathcal{C}}^{(k)}(H,K))  = [H_k : H]^2 l_k$, we get   $ m= [C_k : H]$. Consequently, $$\mathbb{F}_{q}H_{k+1}  \varepsilon_{\mathcal{C}}^{(k+1)}(H,K) \cong  M_{[H_{k+1}: H]}(\mathbb{F}_{q^{l_{k+1}}}),$$ as desired. This completes the proof of the theorem.  \qed \vspace{.5cm} \\ Given a generalized strong Shoda pair $(H, K)$ of $G$, we denote by $R_{(H,K)}$, the subset $R_{n}$ of $\mathcal{C}_{q}(H/K)$ obtained in Proposition \ref{p1}.
\begin{cor}\label{c2}
If $G$ is a generalized strongly monomial group and $S$ is a complete and irredundant set of generalized strong Shoda pairs of $G$, then 	 \begin{description}\item [(i)] $\mathfrak{e}_{\mathcal{C}}(G, H, K)$, $(H, K) \in S$ and $\mathcal{C} \in R_{(H, K)}$ is the complete and irredundant set of primitive central idempotents of $\F G$; 	\item [(ii)]  $\F G \cong  \bigoplus_{(H, K) \in \mathcal{S}} \bigoplus_{\mathcal{C} \in R_{(H, K)}} M_{[G:H]}(\mathbb{F}_{q^{l_{(H,K)}}})$, where $l_{(H, K)}$ is as defined in Theorem \ref{t1}.
		\end{description}  	
\end{cor}
\noindent{\bf Proof.} (i) In view of (\cite{BK2}, Lemma 3), $\mathfrak{e}(G,H,K)$ is a primitive central idempotent of $\mathbb{Q}G$ for every $(H,K) \in S$ and consequently, $\{\mathfrak{e}(G,H,K)~|~(H,K) \in S\}$ is the complete and irredundant set of primitive central idempotents of $\mathbb{Q}G$. Hence, by taking $i=n$ in Proposition \ref{p1}(i), we have \begin{equation}\label{e10}1=\sum_{(H,K) \in S}\overline{\mathfrak{e}(G,H,K)}=\sum_{(H,K) \in S}\sum_{\mathcal{C}\in R_{(H,K)}}\mathfrak{e}_{\mathcal{C}}(G,H,K).\end{equation} Thus, to prove (i) it is enough to show that the primitive central idempotents on the right hand side of eqn (\ref{e10}) are distinct. Suppose $(H,K)$, $(H',K') \in S$, $C \in R_{(H,K)}$ and $C' \in R_{(H',K')}$. If $(H,K)=(H',K')$, then by Proposition \ref{p1}(ii), $\mathfrak{e}_{\mathcal{C}}(G,H,K)$ and $\mathfrak{e}_{\mathcal{C'}}(G,H',K')$ are distinct. Assume that $(H,K)\neq (H',K')$. By Proposition \ref{p1}, $\overline{\mathfrak{e}(G,H,K)}\mathfrak{e}_{\mathcal{C}}(G,H,K)=\mathfrak{e}_{\mathcal{C}}(G,H,K)=\mathfrak{e}_{\mathcal{C}}(G,H,K)\overline{\mathfrak{e}(G,H,K)}$ and \break $
\overline{\mathfrak{e}(G,H',K')}\mathfrak{e}_{\mathcal{C'}}(G,H',K')=\mathfrak{e}_{\mathcal{C'}}(G,H',K')=\mathfrak{e}_{\mathcal{C'}}(G,H',K')\overline{\mathfrak{e}(G,H',K')}$. Since $\mathfrak{e}(G,H,K)\mathfrak{e}(G,H',K')= 0$, it follows immediately that $\mathfrak{e}_{\mathcal{C}}(G,H,K)\mathfrak{e}_{\mathcal{C'}}(G,H',K')= 0$, i.e., $\mathfrak{e}_{\mathcal{C}}(G,H,K)$ and $\mathfrak{e}_{\mathcal{C'}}(G,H',K')$ are distinct. This proves (i). \vspace{.2cm}\\ (ii) From (i), it follows that $$\mathbb{F}_{q}G\cong \bigoplus_{(H,K) \in S}\bigoplus_{\mathcal{C}\in R_{(H,K)}}\mathbb{F}_{q}G\mathfrak{e}_{\mathcal{C}}(G, H, K)$$ and by Theorem \ref{t1}(ii), $\mathbb{F}_{q}G \mathfrak{e}_{\mathcal{C}}(G,H,K)  \cong M_{[G:H]}(\mathbb{F}_{{q}^{l_{(H,K)}}})$. This proves (ii) and completes the proof. \qed
 \section{An algorithm} The theory developed in the previous section provides an algorithmic method to write a complete and irredundant set of primitive central idempotents and the precise Wedderburn decomposition of a generalized strongly monomial group. For this purpose, the following steps are to be followed: \begin{description} \item [Step I]  Write a complete and irredundant set $S$ of generalized strong Shoda pairs of $G$. \item [Step II] For every $(H,K) \in  S$, determine $R_{(H,K)}$ using the following steps: \begin{itemize} \item [(a)] Fix a strong inductive chain $H=H_{0}\leq H_{1}\leq \cdots \leq H_{n}=G$ from $H$ to $G$. \item [(b)] Determine $\mathcal{C}_{q}(H/K)$. \item [(c)] Find a subset $R_{1}$ of $\mathcal{C}_{q}(H/K)$ so that $\{\varepsilon_{\mathcal{C}}^{(0)}(H,K)~|~\mathcal{C} \in R_{1}\}$ is a set of representatives under the action of $\operatorname{Cen}_{H_{1}}(\varepsilon^{(0)}(H,K))$ on $\{\varepsilon_{\mathcal{C}}^{(0)}(H,K)~|~\mathcal{C} \in \mathcal{C}_{q}(H/K) \}$ by conjugation. Note that $R_{1}$ is a set of representatives under the action of $N_{H_{1}}(K)$ on $\mathcal{C}_{q}(H/K)$ by conjugation. \item [(d)] Find a subset $R_{2}$ of $R_{1}$ so that $\{\varepsilon_{\mathcal{C}}^{(1)}(H,K)~|~\mathcal{C} \in R_{2}\}$ is a set of representatives under the action of $\operatorname{Cen}_{H_{2}}(\varepsilon^{(1)}(H,K))$ on $\{\varepsilon_{\mathcal{C}}^{(1)}(H,K)~|~\mathcal{C} \in R_{1} \}$. \item [(e)] Continue the above process for  $n$ steps. Finally, we obtain $R_{n}$ which is the desired $R_{(H,K)}$. \end{itemize} \item [Step III] After finding $R_{(H,K)}$ for every $(H,K) \in S$, one can write primitive central idempotents and Wedderburn decomposition of $\mathbb{F}_{q}G$ using Corollary \ref{c2}. \end{description}  \section{An Example} Let us now illustrate the above algorithm for $G = \operatorname{SmallGroup}(1000,86)$ in GAP library.   Let $\mathbb{F}_{q}$ be a finite field with $q$ coprime to $2$ and $5$ so that $\mathbb{F}_{q}G$ is a semisimple group algebra.  The group $G$ is  generated by $x_{i},~1\leq i\leq 6$, with the following defining relations:
\begin{quote}  $x_{1}^{2}x_{2}^{-1}$=$x_{2}^{2}x_{3}^{-1}$=$x_{4}^{5}$=$x_{3}^{2}$=$x_{5}^{5}$=$x_{6}^{5}$=$1,$\\
               $[x_{2},x_{1}]$=$[x_{3},x_{1}]$=$[x_{3},x_{2}]$=$[x_{6},x_{3}]$=$[x_{6},x_{4}]$=$[x_{6},x_{5}]$=$1,$\\
                $[x_{5},x_{4}]$=$x_{6}$, $[x_{5},x_{1}]$=$x_{4}x_{5},$\\
                $[x_{6},x_{1}]$=$x_{6}^{2}$, $[x_{4},x_{2}]$=$x_{4}x_{6}^{2}$, $[x_{6},x_{2}]$=$x_{6}^{3},$\\
                $[x_{5},x_{2}]$=$x_{5}x_{6}^{2}$, $[x_{5},x_{3}]$=$x_{5}^{3}x_{6}^{2}$, $[x_{4},x_{3}]$=$x_{4}^{3}x_{6}^{2}$, $[x_{4},x_{1}]$=$x_{4}^{2}x_{5}^{3}x_{6}^{4}.$\end{quote} In \cite{BK1}, we have already pointed out that $G$ belongs to the class $\mathcal{C}$. The class $\mathcal{C}$ consists of all finite groups in which each subgroup and each quotient group of subgroup  has a non central abelian normal subgroup.  Also, in Theorem 1 of \cite{BK2}, we have shown  that all the groups in class $\mathcal{C}$ are generalized strongly monomial. Consequently,  $G$ is a generalized strongly monomial group. It is also known that $G$ is not strongly monomial (see \cite{JdR}, Example 3.4.6). Now we proceed to determine the Wedderburn decomposition of $\mathbb{F}_{q}G$.\par The first step is to write a complete and irredundant set of generalized strong Shoda pairs of $G$. In (\cite{BK1}, Section 7), we have seen that  $S$ $=\{ (G,G), \linebreak (G,\langle x_{2},x_{3},x_{4},x_{5},x_{6}\rangle),  (G,\langle x_{3},x_{4},x_{5},x_{6}\rangle), (G,\langle x_{4},x_{5},x_{6}\rangle),  (\langle x_{4},x_{5},x_{6}\rangle,   \langle x_{4},x_{6}\rangle), \linebreak (\langle x_{4},x_{5},x_{6}\rangle,\langle x_{5},x_{6}\rangle),  (\langle x_{4},x_{5},x_{6}\rangle,\langle x_{4}^{-1}x_{5},x_{6}\rangle), (\langle x_{5},x_{6},x_{3}x_{4}^{2}\rangle, \langle x_{3}x_{4}^{2}x_{6}^{3},x_{5}\rangle), \linebreak (\langle x_{5},x_{6},x_{3}x_{4}^{2}\rangle,\langle x_{5}\rangle)\}$ is a complete and irredundant set of Shoda pairs of $G$. Note that all Shoda pairs in $S$ are generalized strong Shoda pairs of $G$. \par The next step is to compute $R_{(H,K)}$  for every $(H,K) \in S$.  Observe that all the Shoda pairs in $S$ except the last two are strong Shoda pairs. Also recall from the algorithm that if $(H, K)$ is a strong Shoda pair of $G$ then $R_{(H, K)}$ is $R_1$, i.e., a set of representatives under the  action of $N_{G}(K)$ on $\mathcal{C}_{q}(H/K)$ by conjugation (by taking strong inductive chain: $H\leq G$).\par The first four Shoda pairs in $S$ are strong Shoda pairs of $G$ and satisfy that  $N_{G}(K)=H$ which gives that the action of $N_{G}(K)$ on $\mathcal{C}_{q}(H/K)$ is trivial. Therefore, if $(H,K) = (G,G),  (G,\langle x_{2},x_{3},x_{4},x_{5},x_{6}\rangle),  (G,\langle x_{3},x_{4},x_{5},x_{6}\rangle)$ or $(G,\langle x_{4},x_{5},x_{6}\rangle)$, then $R_{(H,K)} $ coincides with $ \mathcal{C}_{q}(H/K)$.  Now we compute $ \mathcal{C}_{q}(H/K)$ for these pairs. Clearly, $\mathcal{C}_{q}(G/G)$ contains only the $q$-cyclotomic class of the principal character on $G$ and so there is only one simple component of $\mathbb{F}_{q}G$ corresponding to $(G,G)$ which is isomorphic to $\mathbb{F}_{q}$. For $(H,K)=(G,\langle x_{2},x_{3},x_{4},x_{5},x_{6}\rangle)$, the factor group $H/K$ is cyclic of order $2$ generated by $Kx_{1}$ and hence there is only one $\overline{\mathbb{F}}_{q}$ irreducible character on $H/K$, namely $\sigma: H/K\rightarrow \overline{\mathbb{F}}_{q}$ given by $Kx_{1}\rightarrow \xi_{2}$. Hence $\mathcal{C}_{q}(H/K)$  has only one $q$-cyclotomic class. Consequently, corresponding to $(G,\langle x_{2},x_{3},x_{4},x_{5},x_{6}\rangle)$, there is precisely one simple component of $\mathbb{F}_{q}G$ isomorphic to $\mathbb{F}_{q}(\xi_{2})$, i.e., $\mathbb{F}_{q}$.  Next, for $(H,K) = (G,\langle x_{3},x_{4},x_{5},x_{6}\rangle)$, $H/K$  is cyclic of order $4$ generated by $Kx_{1}$. There are two faithful $\overline{\mathbb{F}}_{q}$ irreducible characters namely $\sigma$ and $\sigma^{3}$ on $H/K$, where $\sigma$ sends $Kx_{1}$ to $\xi_{4}$. Observe that  both $\sigma$ and $\sigma^{3}$ lie in different $q$-cyclotomic class if $q\equiv 1(\operatorname{mod}4)$, and lie in the same $q$-cyclotomic class if $q \not \equiv 1(\operatorname{mod}4)$. Thus, corresponding to $(G,\langle x_{3},x_{4},x_{5},x_{6}\rangle) $, there are two simple components of $\mathbb{F}_{q}G$ (both isomorphic to $\mathbb{F}_{q}(\xi_{4}),$ i.e., $\mathbb{F}_{q}$) if $q\equiv 1\operatorname{mod}4$, and only one simple component (isomorphic to $\mathbb{F}_{q}(\xi_{4})$, i.e., $\mathbb{F}_{q^2}$) if $q \not \equiv 1(\operatorname{mod}4)$. If $(H,K)=(G,\langle x_{4},x_{5},x_{6}\rangle)$, then $H/K$ is cyclic of order $8$ generated by $Kx_{1}$.  Note that $\{ \sigma, \sigma^3, \sigma^5, \sigma^7\}$ is the set of faithful $\overline{\mathbb{F}}_{q}$ irreducible characters on $H/K$, where $\sigma:~H/K\rightarrow \overline{\mathbb{F}}_{q}$ maps $Kx_{1}$ to $\xi_{8}$. Clearly, if $q\equiv 1(\operatorname{mod}8)$, then $\sigma$,$\sigma^3$,$\sigma^5$,$\sigma^7$ lie in distinct $q$-cyclotomic classes. Also, if  $q \not\equiv 1(\operatorname{mod}8)$, then $C_\sigma$, the $q$-cyclotomic class containing $\sigma$, equals $\{\sigma,\sigma^{3}\}$ and $C_{\sigma^5}= \{\sigma^{5},\sigma^{7}\}$. Hence $$\mathcal{C}_{q}(H/K)=\left\{\begin{array}{ll}\{C_\sigma, C_{\sigma^3}, C_{\sigma^5}, C_{\sigma^7}\}, & \hbox{ if~$q\equiv 1(\operatorname{mod}8)$;} \vspace{.2cm}\\\{C_{\sigma}, C_{\sigma^5}\}, & \hbox{if $q \not\equiv 1(\operatorname{mod}8)$.}\end{array}\right.$$ Therefore, corresponding to  $(G,\langle x_{4},x_{5},x_{6}\rangle)$, there are $4$ simple components of $\mathbb{F}_{q}G$ (all isomorphic to $\mathbb{F}_{q}$) if $q\equiv 1(\operatorname{mod}8)$, and two simple components of $\mathbb{F}_{q}G$ (both isomorphic to $\mathbb{F}_{q^{2}}$) if $q \not \equiv 1(\operatorname{mod}8)$. \par The next three Shoda pairs, i.e., $(\langle x_{4},x_{5},x_{6}\rangle, \langle x_{4},x_{6}\rangle),~(\langle x_{4},x_{5},x_{6}\rangle,\langle x_{5},x_{6}\rangle)$ and $(\langle x_{4},x_{5},x_{6}\rangle,\langle x_{4}^{-1}x_{5},x_{6}\rangle)$ are also strong Shoda pairs of $G$. Consider $(H,K)=
(\langle x_{4},x_{5},x_{6}\rangle, \langle x_{4},x_{6}\rangle)$. The factor group $H/K$ is cyclic of order $5$ generated by $Kx_{5}$ and  $$\mathcal{C}_{q}(H/K)=\left\{\begin{array}{lll}\{C_{\sigma},C_{\sigma^{2}},C_{\sigma^{3}},C_{\sigma^{4}}\}, & \hbox{ if~$q\equiv 1(\operatorname{mod}5)$;} \vspace{.2cm}\\\{C_{\sigma}\}, & \hbox{ if~$q\equiv 2{\rm~or~}3(\operatorname{mod}5)$ ;}\vspace{.2cm}\\\{C_{\sigma},C_{\sigma^{2}}\}, & \hbox{if $q \equiv 4(\operatorname{mod}5)$.}\end{array}\right.$$ It can be seen that $N_{G}(K)=\langle x_{2},x_{3},x_{4},x_{5},x_{6}\rangle$. Also, one can check that all the $q$-cyclotomic classes lie in same orbit under the action of $N_{G}(K)$ for all the cases of $q$. Therefore, $R_{(H,K)}=\{C_{\sigma}\}.$ This gives, by Theorem \ref{t1}, that there is only one simple component of $\mathbb{F}_{q}G$ corresponding to the Shoda pair $(\langle x_{4},x_{5},x_{6}\rangle,\langle x_{4},x_{6}\rangle)$ of $G$ and is isomorphic to $M_{8}(\mathbb{F}_{q})$. The similar computations apply to the Shoda pairs $(\langle x_{4},x_{5},x_{6}\rangle,\langle x_{5},x_{6}\rangle)$ and $(\langle x_{4},x_{5},x_{6}\rangle,\langle x_{4}^{-1}x_{5},x_{6}\rangle)$  of $G$ and yield that for each of them there is precisely one simple component (isomorphic to $M_{8}(\mathbb{F}_{q})$) of $\mathbb{F}_{q}G$.
 \par The last two Shoda pairs $ (\langle x_{5},x_{6},x_{3}x_{4}^{2}\rangle,  \langle x_{3}x_{4}^{2}x_{6}^{3},x_{5}\rangle) $ and $ (\langle x_{5},x_{6},x_{3}x_{4}^{2}\rangle,\langle x_{5}\rangle) $ are not strong Shoda pairs of $G$. Let us now compute $R_{(H,K)} $ for these two pairs. Consider $(H, K) = (\langle x_{5},x_{6},x_{3}x_{4}^{2}\rangle,  \langle x_{3}x_{4}^{2}x_{6}^{3},x_{5}\rangle)$. Observe that $H=H_{0}\leq H_{1}\leq H_{2}=G$, where $H_{1}=\langle x_{3},x_{4},x_{5},x_{6}\rangle$, is a strong inductive chain from $H$ to $G$. The factor group $H/K$ is cyclic of order $5$ generated by $Kx_{6}$. Consider $\sigma:~H/K\rightarrow \overline{\mathbb{F}}_{q}$  given by $Kx_{6}\mapsto \xi_{5}$. We have   $$\mathcal{C}_{q}(H/K)=\left\{\begin{array}{lll}\{C_\sigma, C_{\sigma^2}, C_{\sigma^3}, C_{\sigma^4}\}, & \hbox{ if~$q\equiv 1(\operatorname{mod}5)$;} \vspace{.2cm}\\\{C_\sigma \}, & \hbox{ if~$q\equiv 2{\rm~or~}3(\operatorname{mod}5)$;}\vspace{.2cm}\\\{C_\sigma, C_{\sigma^2}\}, & \hbox{if $q \equiv 4(\operatorname{mod}5)$.}\end{array}\right.$$  Since $N_{H_{1}}(K)=H$, the action of $N_{H_{1}}(K)$ on $\mathcal{C}_{q}(H/K)$ is trivial and so $R_{1}=\mathcal{C}_{q}(H/K)$. The next step is to compute $R_2$, which is the desired  $R_{(H,K)}$. Recall from the algorithm that $\{\varepsilon_{\mathcal{C}}^{(1)}(H,K)~|~ \mathcal{C} \in R_{2}\}$ is  a set of representatives under the action of $\operatorname{Cen}_{G}(\varepsilon^{(1)}(H,K))$ on $\{\varepsilon_{\mathcal{C}}^{(1)}(H,K)~|~ \mathcal{C} \in R_{1}\}$. Using Wedderga \cite{wedd}, we can see that $\operatorname{Cen}_{G}(\varepsilon^{(1)}(H,K))=G$. Suppose $q\equiv 1 (\operatorname{mod}5)$. The relations between group elements yield that $\varepsilon_{\mathcal{C}_\sigma}^{(1)}(H,K)^{x_{2}} =  \varepsilon_{\mathcal{C}_{\sigma^4}}^{(1)}(H,K)$. Hence $x_2$  does not centralize  $\varepsilon_{\mathcal{C}_\sigma}^{(1)}(H,K)$. Since $G/H_1$  is cyclic of order $4$ generated by $H_1 x_1$ and $x_1^2 = x_2$, it follows that $\operatorname{Cen}_{G}(\varepsilon_{\mathcal{C_\sigma}}^{(1)}(H,K))=H_{1}$ and so the orbit of $\varepsilon_{\mathcal{C}_\sigma}^{(1)}(H,K)$ has \break $[G:H_1] =4$ elements, namely $\varepsilon_{\mathcal{C}_\sigma}^{(1)}(H,K)$, $\varepsilon_{\mathcal{C}_{\sigma^2}}^{(1)}(H,K)$, $\varepsilon_{\mathcal{C}_{\sigma^3}}^{(1)}(H,K)$ and $\varepsilon_{\mathcal{C}_{\sigma^4}}^{(1)}(H,K)$. Thus $R_{(H, K)} = \{C_\sigma\} .$ Hence, if $q\equiv 1 (\operatorname{mod}5)$, there is only one simple component of $\mathbb{F}_{q}G$, which by Theorem \ref{t1} is isomorphic to $M_{20}(\mathbb{F}_{q})$.
Next, suppose  $q\equiv 2 {\rm~or~}3 (\operatorname{mod}5)$. In this case, $\mathcal{C}_{q}(H/K)$ has only one $q$-cyclotomic class namely $C_\sigma =  \{\sigma, \sigma^{2}, \sigma^{3}, \sigma^{4}\}$. So  $R_{(H,K)}$ equals $ \mathcal{C}_{q}(H/K)$ and in this case also there is one simple component of $\mathbb{F}_{q}G$ isomorphic to $M_{20}(\mathbb{F}_{q})$. Next suppose $q\equiv 4 (\operatorname{mod}5)$.  By Remark \ref{r2},  $\prod_{0\leq i<2}[\operatorname{Cen}_{H_{i+1}}(\varepsilon_{\mathcal{C_\sigma}}^{(i)}(H,K)):H_{i}]$ must divide $2$. In view of Lemma \ref{l3}, $\operatorname{Cen}_{H_1}(\varepsilon_{\mathcal{C_\sigma}}^{(0)}(H,K))$ is a subgroup of $\operatorname{Cen}_{H_1}(\varepsilon^{(0)}(H,K)) = N_{H_{1}}(K) = H$ and therefore, $[\operatorname{Cen}_{H_1}(\varepsilon_{\mathcal{C_\sigma}}^{(0)}(H,K)):H] =1$. Consequently $[\operatorname{Cen}_{G}(\varepsilon_{\mathcal{C_\sigma}}^{(1)}(H,K)):H_{1}] =1$ or $2$. Observe that $x_2$ centralizes $\varepsilon_{\mathcal{C_\sigma}}^{(1)}(H,K) $. This gives that  $[\operatorname{Cen}_{G}(\varepsilon_{\mathcal{C_\sigma}}^{(1)}(H,K)):H_{1}]=2$. Since $[G : H_1]=4$ we have $[G: \operatorname{Cen}_{G}(\varepsilon_{\mathcal{C_\sigma}}^{(1)}(H,K))] =2$. Therefore, the orbit of $\varepsilon_{\mathcal{C_\sigma}}^{(1)}(H,K)$ has two elements, namely $\varepsilon_{\mathcal{C_\sigma}}^{(1)}(H,K)$ and $\varepsilon_{\mathcal{C}_{\sigma^{2}}}^{(1)}(H,K)$. This gives $R_{(H,K)} = \{C_\sigma\} $ and corresponding to $ (\langle x_{5},x_{6},x_{3}x_{4}^{2}\rangle,  \langle x_{3}x_{4}^{2}x_{6}^{3},x_{5}\rangle) $, there is precisely one simple component of $\mathbb{F}_{q}G$ isomorphic to $M_{20}(\mathbb{F}_{q})$. \par Now consider the last Shoda pair $(H,K)=(\langle x_{5},x_{6},x_{3}x_{4}^{2}\rangle, \langle x_{5}\rangle)$.  Observe that  $H=H_{0}\leq H_{1}\leq H_{2}=G$, where $H_{1}=\langle x_{3},x_{4},x_{5},x_{6}\rangle$, is a strong inductive chain. In this case, $H/K$ is cyclic of order $10$ generated by $Kx_{3}x_{4}^{2}$. Suppose $\sigma:~H/K\rightarrow \overline{\mathbb{F}}_{q}$ sends $Kx_{3}x_{4}^{2}$ to $\xi_{10}$. We have $$\mathcal{C}_{q}(H/K)=\left\{\begin{array}{lll}\{C_{\sigma},C_{\sigma^{3}},C_{\sigma^{7}},C_{\sigma^{9}}\}, & \hbox{ if~$q\equiv 1(\operatorname{mod}10)$;} \vspace{.2cm}\\\{C_{\sigma}\}, & \hbox{ if~$q\equiv 3{\rm~or~}7(\operatorname{mod}10)$;}\vspace{.2cm}\\\{C_{\sigma},C_{\sigma^{3}}\}, & \hbox{if $q \equiv 9(\operatorname{mod}10)$.}\end{array}\right.$$
Similar arguments as in the above case  yield that, in all the cases of $q$, there corresponds only one simple component of $\mathbb{F}_{q}G$ isomorphic to $M_{20}(\mathbb{F}_{q})$.   \par  Finally, summing up the above computations, we derive the following Wedderburn decomposition of $\mathbb{F}_{q}G$:

$$\mathbb{F}_{q}G\cong \left\{\begin{array}{lll} \mathbb{F}_{q}^{(8)}\oplus M_{8}(\mathbb{F}_{q})^{(3)}\oplus M_{20}(\mathbb{F}_{q})^{(2)}, & \hbox{ if~$q\equiv 1(\operatorname{mod}8)$;}  \vspace{.3cm}\\
 \mathbb{F}_{q}^{(2)}\oplus \mathbb{F}_{q^{2}}^{(3)}\oplus M_{8}(\mathbb{F}_{q})^{(3)}\oplus M_{20}(\mathbb{F}_{q})^{(2)}, & \hbox{ if~ $q\equiv 3 {\rm~or~}7 (\operatorname{mod}8)$;} \vspace{.3cm}\\ \mathbb{F}_{q}^{(4)}\oplus\mathbb{F}_{q^{2}}^{(2)}\oplus M_{8}(\mathbb{F}_{q})^{(3)}\oplus M_{20}(\mathbb{F}_{q})^{(2)}, & \hbox{ if~
 $q\equiv 5 (\operatorname{mod}8).$}\end{array}\right.$$ \\
Here, $R^{(n)}$ denotes  the sum of $n$ copies of $R$.

 \bibliographystyle{amsplain}
\bibliography{Reference}

\end{document}